\documentclass[12pt]{amsart}

\usepackage[dvips]{graphicx}
\usepackage{amsfonts,amssymb,amscd,bm}
\newcommand{\R}{{\mathbb R}}

\newcommand{\Z}{{\mathbb Z}}

\newcommand{\e}{{\varepsilon}}

\newcommand{\al}{{\alpha}}
\newcommand{\be}{{\beta}}
\newcommand{\de}{{\delta}}

\newcommand{\Ga}{{\Gamma}}

\newcommand{\ti}{\tilde}
\newcommand{\id}{{\mathrm{id}}}

\newcommand{\vl}{{\; | \;}}

\newcommand{\bu}{{\bullet}}

\newcommand{\ab}[1]{ {\langle #1\rangle} }

\newcommand{\cusp}{ \begin{picture}(10,10)(0,0)
  \put(1,0){\oval(8,12)[tr]}
  \put(9,0){\oval(8,12)[tl]}
 \end{picture} }

\newcommand{\tang}{ \begin{picture}(10,10)(0,2)
  \put(1,5){\oval(8,8)[r]}
  \put(9,5){\oval(8,8)[l]}
 \end{picture} }

\newcommand{\trip}{ \begin{picture}(11,10)(0,2)
  \put(0,0){\line(1,1){10}}
  \put(5,0){\line(0,1){10}}
  \put(0,10){\line(1,-1){10}}
 \end{picture} }

\theoremstyle{definition}

\title[Gauss paragraphs of classical links]
{Gauss paragraphs of classical links and\\
 a characterization of virtual link groups}

\author[kurlin]{V.~Kurlin$^*$}
\address{Department of Mathematical Sciences,   
 The University of Liverpool,
 L69 7ZL, Liverpool, United Kingdom}
\email{ kurlin@liv.ac.uk }

\date{First version: October 30, 2006. This version: February 18, 2007. }

\begin{document}
\vspace*{-5mm}

\maketitle

\begin{abstract}
A classical link in 3-space can be represented 
 by a Gauss paragraph encoding a link diagram
 in a combinatorial way.
A Gauss paragraph may code not a classical link diagram,
 but a diagram with virtual crossings.
We present a criterion and a linear algorithm detecting 
 whether a Gauss paragraph encodes a classical link.
We describe Wirtinger presentations
 realizable by virtual link groups.
\end{abstract}


\section{Introduction}


\subsection{Brief summary}
\noindent
\smallskip

This is a research on the interface between knot theory, 
 group theory and combinatorics.
Briefly, a classical knot is a closed loop without self-intersections,
 considered up to a smooth deformation of the ambient 3-sphere.
A classical link is a union of several disjoint closed loops.
\smallskip

A plane diagram of a classical link can be encoded 
 combinatorially by a Gauss diagram, a union of circles with arrows.
Each arrow connects two points that map to a crossing of the plane diagram.
Signs and orientations of arrows allow us 
 to specify the overcrossing information.
\smallskip

Plane diagrams of links are defined up to Reidemeister moves
 converting into moves on Gauss paragraphs,
 word codes of Gauss diagrams.
A Gauss paragraph may code not a classical link diagram, 
 but a diagram with virtual crossings
 without specified overcrossing information.
\smallskip

Virtual links generalize classical ones and are defined via 
 plane diagrams with classical and virtual crossings up to natural moves.
Geometrically, a virtual link can be considered as 
 a union of linked closed loops in a thickened surface.
Virtual crossings invisible on the surface 
 appear in a plane diagram when the surface is projected to a plane.
\smallskip

By definition a Wirtinger group has a Wirtinger presentation,
 where each relation says that two generators are conjugate.
The fundamental groups of classical knot complements have
 Wirtinger presentations with the additional restriction 
 that all generators are conjugate \cite{CF}.
\smallskip

The notion of the fundamental group for classical knot complements
 extends to a larger class of virtual knots.
The resulting Wirtinger groups possess unusual properties,
 e.g. have non-trivial second homology groups. 
The groups of virtual knots can be characterized as groups with
 Wirtinger presentations, where all generators are 
 conjugate and their total number either equals the number of 
 defining relations or exceeds it by 1 \cite{Kim}.


\subsection{Classical and virtual links}
\noindent
\smallskip

Here we introduce basic notions of the classical knot theory
 and its generalization to virtual links proposed by L.~Kauffman in \cite{Kau}.
\medskip

\noindent
{\bf Definition 1.1.}
A \emph{classical knot} is the image of a smooth embedding $S^1\to S^3$,
 i.e. a closed loop without self-intersections, see Fig.~1.
A \emph{classical link} is a smooth embedding of 
 several disjoint circles.
\medskip

We will consider oriented links with unordered components.
Usually links are studied up to
 isotopy that is a smooth deformation of $S^3$,
 see Definition~2.1.
Classical links are represented by plane diagrams 
 defined up to Reidemesiter moves, see Proposition~2.2.
The classification problem of knots up to isotopy 
 has an algorithmic solution whose complexity is 
 highly exponential in the number of crossings \cite{Mat}.
A polynomial algorithm is not known even 
 for recognizing the \emph{unknot}, a round circle.

\begin{figure}[!h]
\includegraphics[scale=1.0]{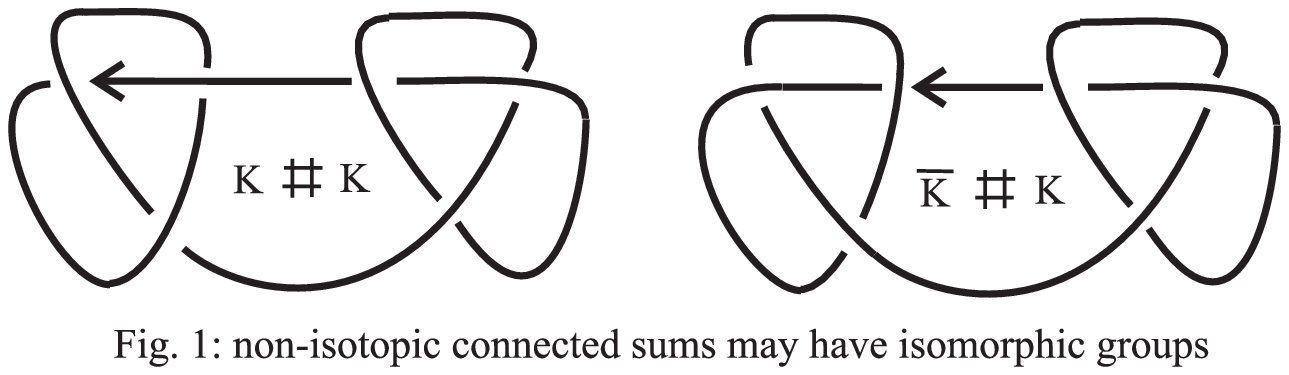}
\end{figure}

The fundamental groups of knot complements are powerful invariants
 distinguishing all \emph{prime} knots that are not connected sums
 of non-trivial knots, see Fig 1.
The fundamental group with an additional peripheral structure
 is a complete knot invariant \cite{Wal}.
But this algebraic classification does not provide 
 an effective algorithm for detecting knots.
\medskip

\noindent
{\bf Definition 1.2.} 
A \emph{virtual} link is a smooth immersion of several oriented circles
 into the plane, that is a smooth embedding outside finitely many 
 double transversal intersections of 2 types, see Fig.~2:
\smallskip

\noindent
$\bu$ \emph{classical} crossings, where one arc overcrosses the other;

\noindent
$\bu$ \emph{virtual} crossings, where 
 the intersecting arcs are not distinguishable.
\smallskip

\noindent
The \emph{plane} diagram of a virtual link is its image in the plane
 with the specified overcrossing information at classical crossings only.

\begin{figure}[!h]
\includegraphics[scale=1.0]{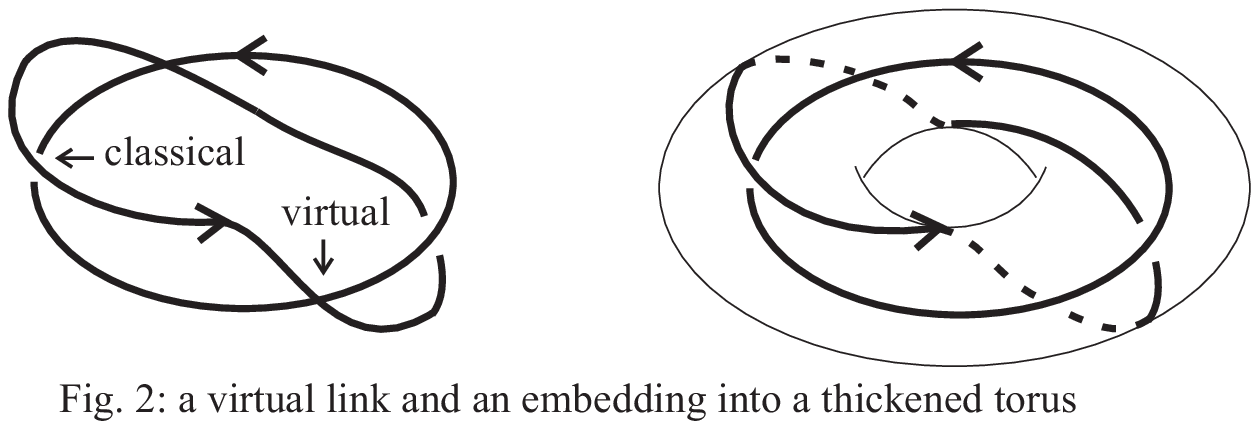}
\end{figure}

Virtual knot theory is motivated by studying knots in 
 thickened surfaces $S_g\times\R$, see \cite{JKS}.
Virtual crossings appear under a projection $S_g\to\R^2$.
Given a diagram of a virtual knot, 
 the genus of $S_g$ is unknown.
The equivalence relation for virtual knots is 
 introduced via moves in Definition~4.1.


\subsection{Gauss paragraphs of links}
\noindent
\smallskip

Firstly we define an abstract Gauss paragraph
 as a collection of words.
After we encode the plane diagram of a classical link 
 by a Gauss paragraph.
Our first contribution is Algorithm~1.4 computing
 the least genus of a thickened surface containing 
 a link encoded by a given Gauss paragraph.
We also prove a criterion for the planarity of a Gauss paragraph
 in Theorem~3.6.
Our second contribution is a characterization of Wirtinger presentations
 realizable as the groups of virtual links, see Theorem~4.8.
\medskip

\noindent
{\bf Definition 1.3.}
Given $n\geq 1$ fix the alphabet $\{i,\, i^+,\, i^- \vl i=1,\dots,n\}$.
A \emph{Gauss paragraph} is an unordered collection $\{u_1,\dots,u_k\}$ 
 of words such that every word is of even length and,
 for each $i\in\{1,\dots,n\}$, exactly one letter $i$ 
 and either $i^+$ or $i^-$ belong to 
 the union $u_1\cup\dots\cup u_k$. 

\begin{figure}[!h]
\includegraphics[scale=1.0]{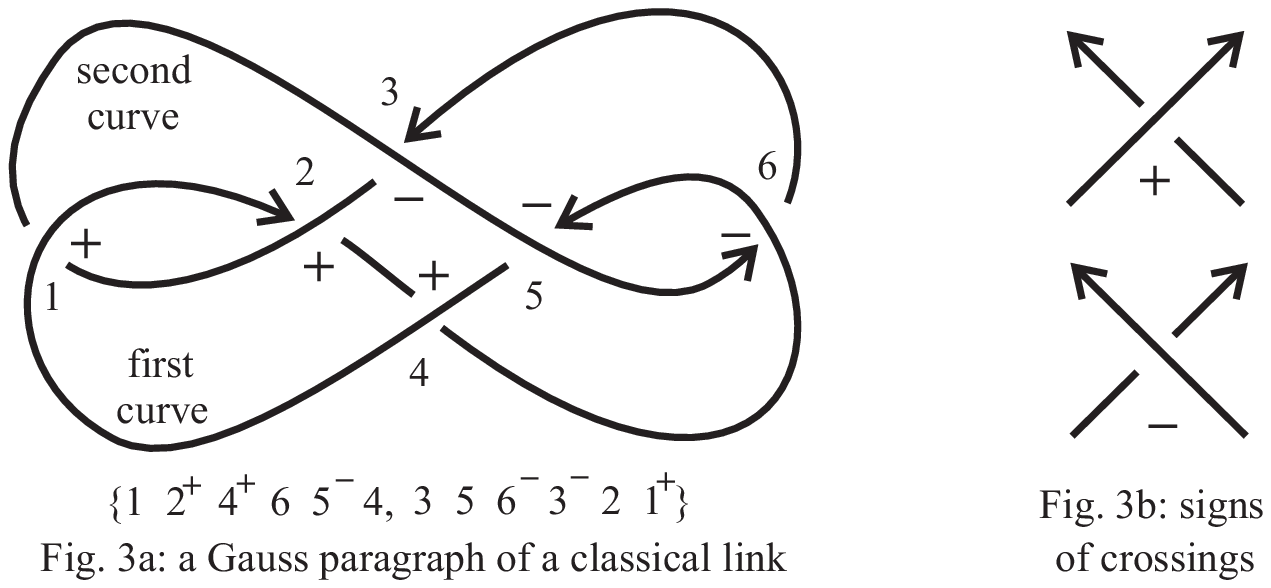}
\end{figure}

Enumerate all crossings in the diagram of 
 a classical oriented link in any order, see Fig.~3a.
To each curve in the diagram associate a word as follows.
Write down the numbers of crossings according to their order
 in the curve.
If we pass through an undercrossing $i$, we add 
 the superscript $\e$ to the letter $i$, where
 $\e$ is the \emph{sign} of the crossing, see Fig~3b.
The resulting collection of words is a Gauss paragraph
 of the plane diagram.
The construction is the same for diagrams
 in oriented surfaces.
Any link diagram can be encoded by several Gauss paragraphs since
 we can renumber the crossings and permute the components.
\smallskip

Any Gauss paragraph gives rise to a link embedded into 
 a thickened surface since any virtual crossing disappears 
 after adding a 1-handle, see Fig.~2.
In the least genus surface the embedded link is unique 
 up to homeomorphisms of the surface.
The virtual knot theory coincides with the theory of links
 in thickened surfaces up to addition and subtraction of 1-handles
 outside the link, see \cite{Kau}.
\medskip

\noindent
{\bf Algorithm~1.4}.
There is an algorithm of linear complexity $Cn$
 to compute the least genus of an oriented surface
 containing a link diagram encoded by a given Gauss paragraph
 consisting of $2n$ letters. 
\medskip

\noindent
{\bf Acknowldegement.}
The author thanks L.~Kauffman, H.~Morton, C.~Souli\'e for useful discussions
 and D.~Elton for his manuscript \cite{CE}.
The author was supported by
 Marie Curie Fellowship 007477.


\section{Plane diagrams and Gauss diagrams}


\subsection{An isotopy of knots and Reidemeister moves}
\noindent
\smallskip

\noindent
{\bf Definition 2.1.}
Links $K,L\subset S^3$ are called \emph{isotopic} if 
 there is a continuous family of homeomorphisms $F_t:S^3\to S^3$,
 $t\in[0,1]$, an \emph{isotopy}, such that $F_0=\id$ and $F_1(K)=L$.
\medskip

Plane diagrams of links are supposed to be \emph{regular}, i.e. 
 they have finitely many double transversal intersections, 
 \emph{crossings}.
An \emph{isotopy} of diagrams is a smooth family of 
 regular diagrams.
In generic isotopies of diagrams only the following codimension~1 singularities
 can occur: an ordinary cusp $\cusp$, a simple tangency $\tang$ 
 and a transversal triple intersection $\trip$.
The Reidemeister theorem below is an application of singularity theory
 stating that any isotopy of knots can be made transversal to 
 the three codimension~1 discriminants in the space of all smooth links.
\medskip

\noindent
{\bf Proposition 2.2.} 
(Reidemeister theorem)
Two regular diagrams represent isotopic links if and only if
 they can be connected by isotopies of diagrams and
 finitely many \emph{Reidemeister} moves in Fig.~4.
\qed

\begin{figure}[!h]
\includegraphics[scale=1.0]{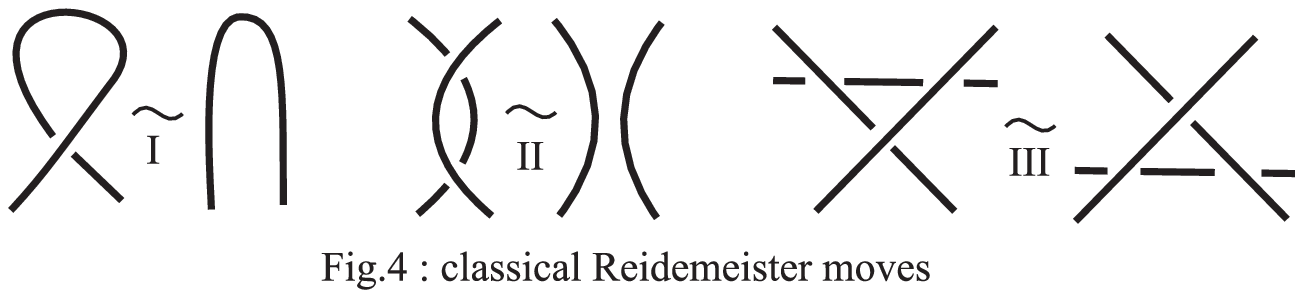}
\end{figure}

In the above Reidemeister theorem we should consider 
 all symmetric images of the moves in Fig.~4 
 for all possible orientations.


\subsection{Virtual links via plane diagrams}
\noindent
\smallskip

Firstly we define the equivalence relation for virtual links
 as in \cite{Kau}.
\medskip

\noindent
{\bf Definition 2.3.} 
Two virtual links are \emph{equivalent} if their 
 plane diagrams can be connected by a finite sequence 
 of the following moves
\smallskip

\noindent
$\bu$ 
 three \emph{classical} Reidemeister moves in Fig.~4;

\noindent
$\bu$ 
 three \emph{virtual} Reidemeister moves, 
 where all crossings are virtual;

\noindent
$\bu$ 
 the \emph{mixed} move in Fig.~5a with 
 classical and virtual crossings.
\medskip

\begin{figure}[!h]
\includegraphics[scale=1.0]{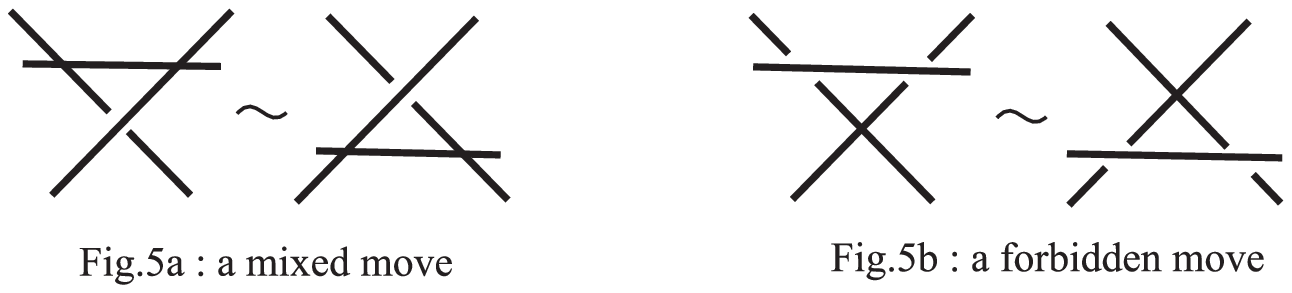}
\end{figure}

The move in Fig.~5b and its symmetric images are \emph{forbidden}, 
 because any virtual link can be transformed to the \emph{unlink}, 
 the disjoint union of circles, through the moves of 
 Definition~2.3 and Fig.~5b, see \cite{GPV}.
The moves of Definition~2.3 generate formally a new equivalence 
 relation on classical link diagrams without virtual crossings,
 but the resulting theory coincides with the theory of classical links, see \cite{GPV}.


\subsection{Virtual links via Gauss diagrams}
\noindent
\smallskip

Now we shall look at virtual links from a purely combinatorial point of view.
Virtual links can be represented by Gauss diagrams.
\medskip

\noindent
{\bf Definition 2.4.} 
A \emph{Gauss} diagram consists of (see Fig.~6):

\noindent
$\bu$ 
 a union of several oriented circles;

\noindent
$\bu$ 
 arrows connecting points on these circles;

\noindent
$\bu$ 
 a sign $+$ or $-$ associated to each arrow.
\smallskip

\noindent
Two Gauss diagrams are considered up to orientation preserving
 diffeormorphism respecting the arrows and their signs.
\medskip

\begin{figure}[!h]
\includegraphics[scale=1.0]{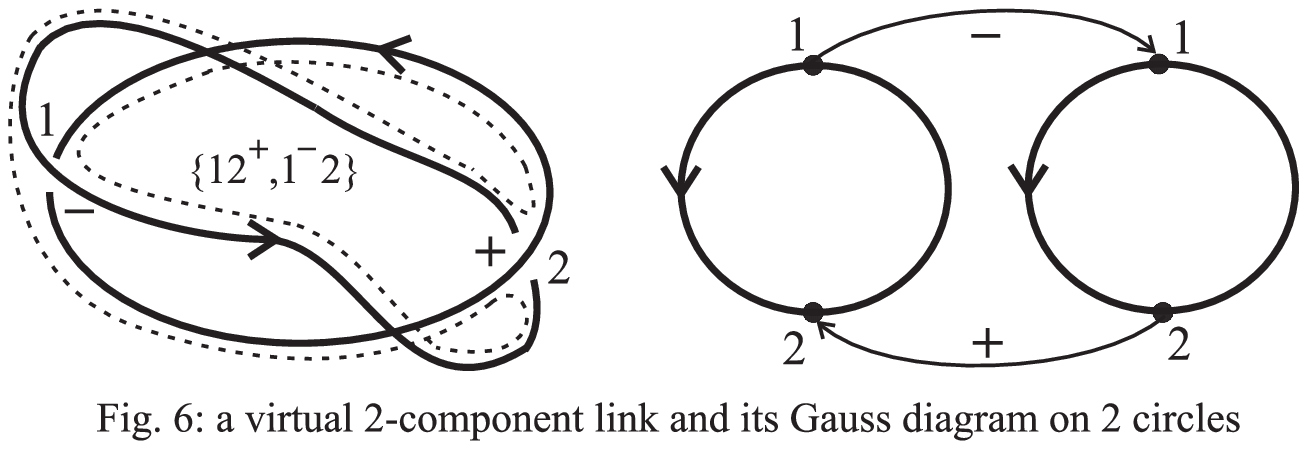}
\end{figure}

Any plane diagram $D$ of a $k$-component virtual link $K$ 
 gives rise to the following Gauss diagram on $k$ circles.
The components of $K$ are embeddings of $k$ circles.
Two points on these circles are connected by an arrow if and only if
 they map to a classical crossing in $D$.
The arrow is oriented from the upper branch to the lower one in $D$.
Each arrow is equipped with the sign of the crossing, see Figs.~3b and 6.
\smallskip

The moves in Fig.~4 generate 
 the transformations of Gauss diagrams in Fig.~7.
The round arcs there may belong to different circles.
The Reidemeister move III in Fig.~4 gives rise to 
 eight transformations with specified orientations.
But the six transformations in Fig.~7 are sufficient to realize
 all Reidemeister moves in Fig.~4 \cite{Ost}.
\smallskip

\begin{figure}[!h]
\includegraphics[scale=1.0]{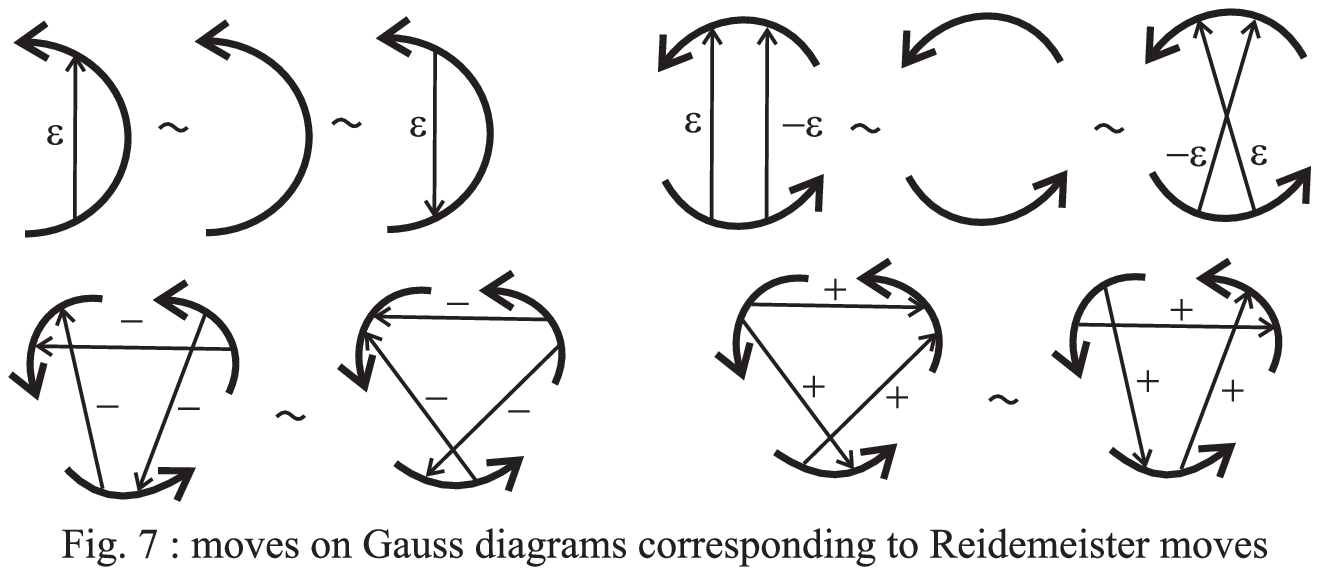}
\end{figure}

The transformations in Fig.~7 can be realized by 
 the moves of Definition~2.3 on plane diagrams.
The virtual Reidemeister moves and the mixed move in Fig.~5a
 were designed so that any arc involving virtual crossings only
 can be replaced in a given diagram by an arbitrary arc having
 the same endpoints and intersecting other arcs in virtual crossings.
So plane diagrams of virtual links up to the moves of Definition~2.3
 are equivalent to Gauss diagrams up to the transformations in Fig.~7.
\smallskip

Any Gauss diagram on $k$ circles can be converted into 
 a Gauss paragraph as follows.
Number the arrows of the Gauss diagram by $1,\dots,n$ in any order.
Label the tail and head of the $i$th arrow with a sign $\e$ 
 by $i$ and $i^{\e}$, respectively.
By reading the labels counterclockwisely
 we obtain words $u_1,\dots,u_k$ considered cyclically, see Fig.~6.


\section{The planarity of Gauss paragraphs}


\subsection{A linear algorithm for the planarity of Gauss paragraphs}
\noindent
\smallskip

Any Gauss paragraph can be realized by a diagram with 
 classical crossing on a suitable oriented surface.
So any Gauss paragraph of $k$ words encodes 
 an embedding of $k$ circles into a thickened surface.
We shall compute the least genus of such a surface
 in linear time with respect to the length of the Gauss paragraph.
If the genus is 0, Algorithm~1.4
 determines whether a Gauss paragraph encodes a classical link.
\medskip

\noindent
{\bf Definition 3.1.}
Let $\{u_1,\dots,u_k\}$ be a Gauss paragraph consisting of $2n$ letters.
We shall construct its \emph{Carter surface} $M\{u_1,\dots,u_k\}$ as 
 a combinatorial cell complex \cite{Car}.
Take $n$ vertices labelled by $1,\dots,n$.
\smallskip

We connect vertices $i,j$ by an edge with a mark $(a,b)$ or $(b,a)$, 
 if one of the cyclic words $u_1,\dots,u_k$ contains
 the ordered pair $ab$ or $ba$ of successive letters, respectively, 
 for some $a\in \{i,i^+,i^-\}$, $b\in \{j,j^+,j^-\}$.
For instance, the Gauss paragraph $\{12^+, 1^-2\}$ in Fig.~6 
 generates the graph with 2 vertices connected by 4 edges
 $(12^+),(2^+1),(1^-2),(21^-)$.
\smallskip

Let us select unoriented cycles in the resulting graph.
Travelling along an edge $(a,b)$, encode the direction of our path
 by $(a,b)_+$ if the letter $a$ preceeds $b$ in the Gauss paragraph,
 otherwise by $(a,b)_-$.
After passing an edge, we choose the next one 
 by the following local rules (Fig.~8):
\smallskip

\noindent
$(a,i)_+ \to (i^{\de},b)_{\de}$, 
$(a,i)_- \to (i^{\de},b)_{-\de}$,
$(a,i^+)_{\e} \to (i,b)_{-\e}$, 
$(a,i^-)_{\e} \to (i,b)_{\e}$ 
\smallskip

\noindent
for a unique possible choice of $\de=\pm$ and both values of $\e=\pm$.
\medskip

\begin{figure}[!h]
\includegraphics[scale=1.0]{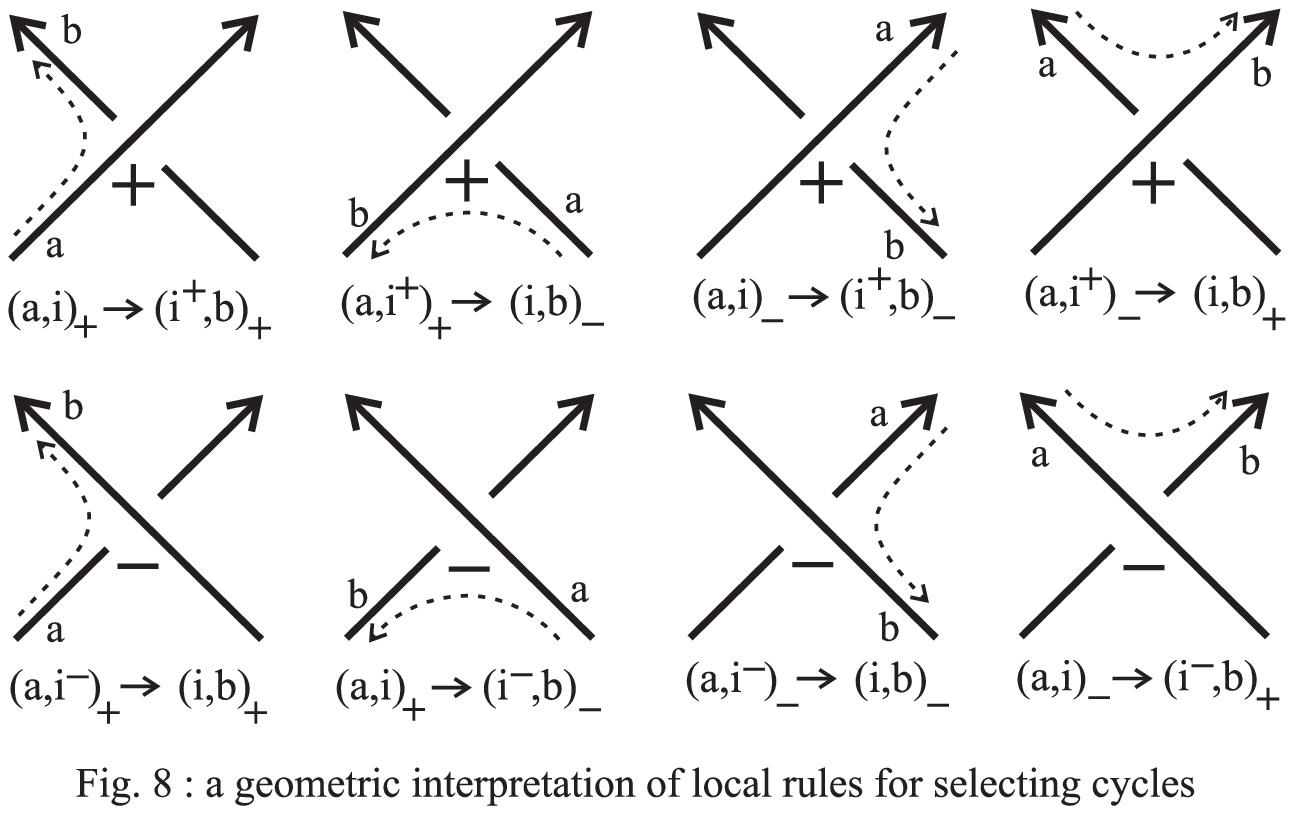}
\end{figure}

A geometric interpretation of the above pure combinatorial rules 
 is in Fig.~8.
If the given Gauss paragraph encodes a diagram on an oriented surface,
 then the orders of letters in the cyclic words correspond 
 to the orientations of the components.
Geometrically, the rules say that we always turn left.
After passing through each edge once in both directions, 
 we have selected all cycles defining faces in the Carter surface.
\smallskip

The Gauss paragraph $\{12^+, 1^-2\}$ in Fig.~6 
 leads to 2 faces bounded by the cycles 
 $(12^+)_+(21^-)_-(12^+)_-(21^-)_+$ and $(2^+1)_+(1^-2)_-(2^+1)_-(1^-2)_+$.
The first cycle is shown by the dashed closed curve in Fig.~6.
So the Carter surface $M\{12^+,1^-2\}$ has Euler characteristic 0 and is a torus. 
\medskip

\noindent
{\bf Lemma 3.2.}
Given a Gauss paragraph the Carter surface
 has the least genus among all oriented surfaces
 containing a classical link diagram encoded by the Gauss paragraph.
A Gauss paragraph encodes a classical link 
 if and only if the Carter surface is a sphere.
\begin{proof}
The Carter surface of the given Gauss paragraph contains a desired link diagram 
 due to the geometric interpretation of Definition~3.1 in Fig.~8. 
Suppose that there is another surface $S_g$ containing a link diagram
 encoded by the Gauss paragraph.
The underlying graph of the diagram with classical crossings only
 splits $S_g$ into several connected pieces.
The genus $g$ is minimal, if all the pieces are disks
 as in Definition~3.1.
\end{proof}
\medskip

\noindent
\emph{Proof of Algorithm~1.4.}
Assume that a Gauss paragraph $\{u_1,\dots,u_k\}$ 
 encodes a connected Gauss diagram.
Denote by $|u_1|,\dots,|u_k|$ the lengths of the words.
By the rules of Definition~3.1 
the genus $g$ the Carter surface $M\{u_1,\dots,u_k\}$ can be computed
 in linear time with respect to $2n=|u_1|+\dots+|u_k|$ via 
 $\chi(M\{u_1,\dots,u_k\})=n-2\sum_{i=1}^k|u_i|+\#$(faces)$=2-2g$.
By Lemma~3.2 the resulting genus $g$ is minimal.
If $g=0$, the Gauss paragraph represents a classical link.
If the Gauss paragraph splits into disjoint subparagraphs,
 run the algorithm for each of them.
\qed


\subsection{A criterion for the planarity of Gauss paragraphs}
\noindent
\smallskip

To formulate the criterion we first define abstract Gauss codes
 and then associate a Gauss code to each Gauss paragraph.
\medskip

\noindent
{\bf Definition 3.3.}
A \emph{Gauss code} $W$ is a permutation of the symbols
 $1^{+1},1^{-1},\dots,n^{+1},n^{-1}$, considered as a cyclic word.
Let $S_i$ be the subword of $W$ between $i^{+1}$ and $i^{-1}$,
 not including these symbols.
For $i=1,\dots,n$, let $\al_i(W)$ be the sum of the superscripts
 of the symbols of $S_i$.
\smallskip

Denote by $S_i^{-1}$ the set of the symbols of $S_i$, 
 where all superscripts are reversed.
Put $\bar S_i=S_i\cup\{i^{+1},i^{-1}\}$.
For $i,j\in\{1,\dots,n\}$, let $\be_{ij}(W)$ be 
 the sum of the superscripts
 of the symbols of $\bar S_i\cap S_j$.
\medskip

A generic immersion of an oriented circle into the plane 
 can be represented by a Gauss code as follows.
Attach indices $1,\dots,n$ to the crossings in any order, 
 go along the curve and write down the corresponding indices.
We add the superscript $+1$ to an index $i$ if 
 the crossing branch at the $i$th crossing goes 
 from left to right, otherwise $-1$, see Fig.~10.
\medskip

\noindent
{\bf Proposition~3.4} \cite{CE}
A Gauss code $W$ encodes a planar oriented closed curve
 if and only if $\al_i(W)=\be_{ij}(W)=0$ for all $i,j\in\{1,\dots,n\}$.

\begin{figure}[!h]
\includegraphics[scale=1.0]{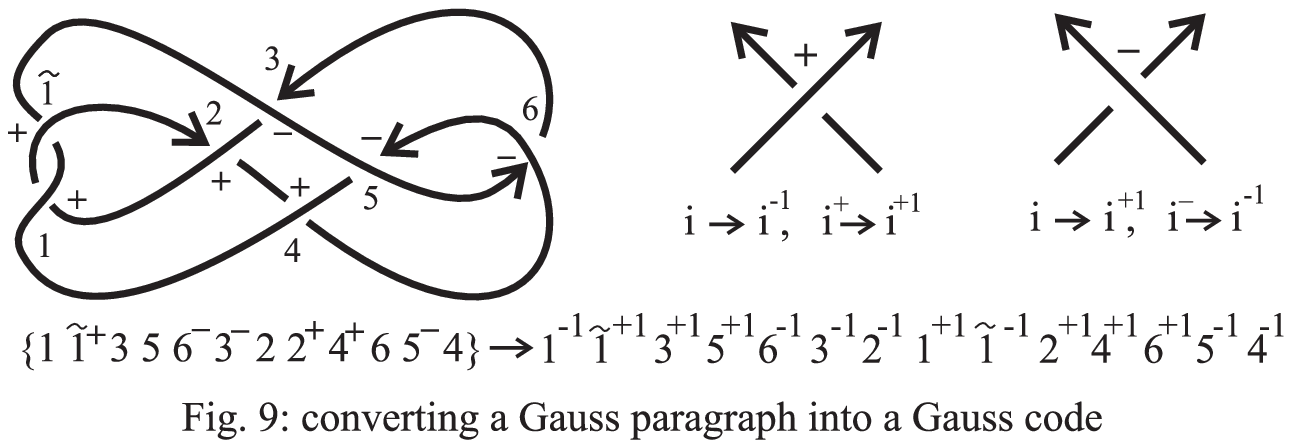}
\end{figure}

For simplicity we assume that a Gauss paragraph can not be split into
 a disjoint union of non-empty subparagraphs.
Construction~3.5 and Theorem~3.6 extends trivially 
 to the case of splittable paragraphs.
\medskip

\noindent
{\bf Construction 3.5.}
To a Gauss paragraph $\{u_1,\dots,u_k\}$ consisting of $2n$ letters
 we associate a Gauss code $W$ of length $2n+2k-2$.
Take two words, say $u_1,u_2$, such that $u_1$ contains
 a letter $i$ and $u_2$ has $i^+$ or $i^-$.
\smallskip

Rewrite $u_2$ cyclically in such a way that 
 $i^+$ or $i^-$ is the last letter of $u_2$, say $i^+$.
Double these leters by the rule:
 $i\mapsto i \tilde i^+$ and $i^+\mapsto i^+ \tilde i$.
If the sign of the $i$th crossing is negative, use the rules
 $i\mapsto i {\tilde i}^-$ and $i^-\mapsto i^- \tilde i$.
Now insert the word $u_2$ ending by $i^+ \tilde i$
 right after $i \tilde i^+$ into the word $u_1$.
We get a new word of length $|u_1|+|u_2|+2$.
The Gauss paragraph of 2 words in Fig.~3 becomes
 a single word in Fig.~9 at the left.
\smallskip

Continue uniting words until we get a single word of length $2n+2k-2$.
Finally replace letters as follows:
 $i^+\to i^{+1}$, $i^-\to i^{-1}$, $i\to i^{-\e}$, where
 $\e$ is the superscipt of the other letter with the same index $i$,
 see Fig.~9.
\medskip

A Gauss code $W$ associated to a given Gauss paragraph is not 
 uniquely defined, but the following criterion works always.
\medskip

\noindent
{\bf Theorem 3.6.}
A Gauss paragraph $\{u_1,\dots,u_k\}$ consisting of $2n$ letters 
 encodes the plane diagram of a classical $k$-component link
 if and only if the invariants $\al_i$ and $\be_{ij}$, $i,j=1,\dots,n$, 
 vanish for any Gauss code $W$ associated to the given Gauss paragraph 
 in Construction~3.5.
\begin{proof}
We prove that the given Gauss paragraph is planar if and only if
 the Gauss code obtained via Construction~3.5 is planar.
\smallskip

Let $D_U$ and $D_W$ be the plane diagrams encoded by 
 the Gauss paragraph and Gauss code, respectively.
If $D_U$ has only classical crossings then $D_W$ is obtained from $D_U$
 by doubling $k-1$ crossings, where different components intersect.
Conversely, if $D_W$ has only classical crossings then
 $D_U$ is obtained from $D_W$ by compressing $k-1$ pairs of crossings
 into $k-1$ single crossings.
It remains to apply Proposition~3.4.
\end{proof}

\begin{figure}[!h]
\includegraphics[scale=1.0]{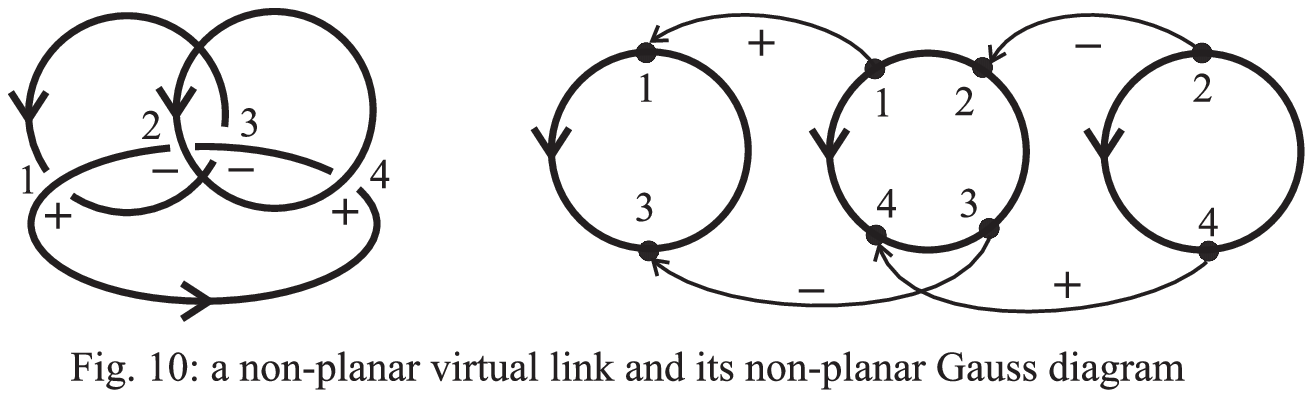}
\end{figure}

More complicated criteria for the planarity of Gauss paragraphs
 were obtained in \cite{Sch, Sou} via 
 associated graphs and word partitions.
Fig.~10 contains a virtual analogue of the Brunnian link.
Any pair of the 3 components is planar, i.e. the corresponding Gauss diagram
 on 2 circles encodes a plane diagram with classical crossings only.
Theorem~3.6 implies immediately that the whole virtual link 
 on 3 components is not planar since after doubling 
 crossings 1 and 4 the resulting Gauss code 
 $1^{-1}4^{+1}\ti 4^{-1}2^{+1}4^{-1}\ti 4^{+1}3^{+1}2^{-1}\ti 1^{-1}1^{+1}3^{-1}\ti 1^{+1}$
 has $\al_2\neq 0$, see Definition~3.3.


\section{The realizability of Wirtinger groups}


\subsection{Virtual link groups}
\noindent
\smallskip

Firstly we introduce abstract Wirtinger presentations.
Then to each virtual link we associate the group
 with a Wirtinger presentation.
\medskip

\noindent
{\bf Definition 4.1.} 
A \emph{Wirtinger} group is a group with a \emph{Wirtinger} presentation
 $\Pi=\ab{m_1,\dots,m_n\vl r_1,\dots,r_s}$, where $n\geq s$ 
 and each relator $r_q$ has the form $m_i=w_q^{-1}m_jw_q$
 for some $i,j\in\{1,\dots,n\}$ and words $w_q$ in 
 the generators $m_1,\dots,m_n$, $q=1,\dots,s$.
The \emph{weight} of a group is the minimal number $k$ 
 of elements whose conjugates generate the group.
\medskip

The fundamental group $\pi(L)=\pi_1(S^3-L)$ 
 of the complement to a classical $k$-component link 
 always has a Wirtinger presentation \cite{CF}.
The generators of $\pi(L)$ split into $k$ conjugate classes, 
 i.e. the link group $\pi(L)$ is a Wirtinger group of weight $k$.
So the \emph{abelianization} $\pi(L)/[\pi(L),\pi(L)]$ is 
 isomorphic to $\Z^k$ for any $k$-component link $L$.
We give a general construction for the group of a virtual link
 from a Gauss diagram.
\medskip

\noindent
{\bf Definition 4.2.} 
Given a Gauss diagram $G$ with $k$ circles and $n$ arrows,
 a group $\Pi$ with $n$ generators and $n$ relations
 will be constructed.
After cutting the circles at each arrowhead and 
 forgetting all arrowtails, $G$ splits into $n$ arcs.
To each arc associate a generator of $\Pi$, see Fig.~11.
\smallskip

Then every arrow gives rise to a defining relation in $\Pi$.
Let $\e$ be the sign of an arrow in $G$, the tail of the arrow lie
 on the arc denoted by $a$.
If the head of the arrow is between the arcs $b$ and $c$
 in the counterclockwise direction then 
 the associated relation is $c=a^{-\e}ba^{\e}$, see Figs.~6, 11.

\begin{figure}[!h]
\includegraphics[scale=1.0]{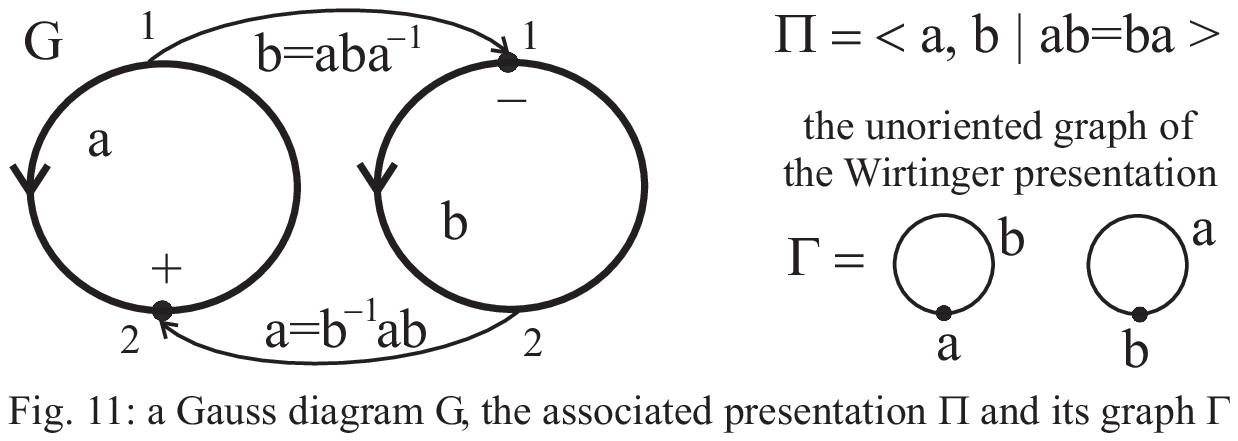}
\end{figure}

By Definitions~4.1 and 4.2 the presentation associated to 
 a Gauss diagram is a Wirtinger one.
The \emph{group} of a virtual link is the group $\Pi$ of 
 a Gauss diagram representing the virtual link, see Fig.~6.
Since generators of $\Pi$ corresponding to adjacent arcs
 are conjugate, Lemma~4.3 follows.
\medskip

\noindent
{\bf Lemma 4.3.}
Let $\Pi$ be the group of a Gauss diagram $C$ on $k$ circles.
Then all generators corresponding to the arcs of same circle in $C$
 are conjugate to each other. 
Hence $\Pi$ is normally generated by $k$ elements,
 the abelianization $\Pi/[\Pi,\Pi]$ is isomorphic to $\Z^k$.
\qed


\subsection{The graph of a Wirtinger group}
\noindent
\smallskip

All generators of a Wirtinger group split into
 classes of congugate ones.
More explicitly the structure of a Wirtinger group
 can be described by a finite graph $\Ga$ associated 
 to its presentation $\Pi$.
\medskip

\noindent
{\bf Definition 4.4.} 
A Wirtinger presentation
 $\Pi=\ab{m_1,\dots,m_n\vl r_1,\dots,r_s}$, where
 $r_q=m_i^{-1} w_q^{-1} m_j w_q$,
 defines the \emph{graph} $\Ga$ as follows.
The vertices of $\Ga$ are in a 1-1 correspondence
 with the generators $m_1,\dots,m_n$.
Vertices $m_i,m_j$ are connected by an edge marked by $w_q$
 if and only if $\Pi$ has a relator $r_q=m_i^{-1} w_q^{-1} m_j w_q$
 for a word $w_q$ in $m_1,\dots,m_n$, see Fig.~11.
\medskip

The graph $\Ga$ is unoriented and can be disconnected as
 the example in Fig.~11.
The \emph{Euler} characteristic of a connected graph is
 the number of vertices minus the number of edges,
 so it is always not bigger than 1.
\medskip
 
\noindent
{\bf Proposition 4.5.}
Let $\Ga$ be the graph of a Wirtinger presentation $\Pi$.
The connected components of $\Ga$ are in a 1-1 correspondence
 with the classes of conjugate generators of $\Pi$.
If $\Pi$ is a virtual link group then each connected component of $\Ga$ 
 has Euler characteristic 0 or 1.
\begin{proof}
For each relator $r_q=m_i^{-1} w_q^{-1} m_j w_q$ of $\Pi$,
 mark the corresponding edge of $\Ga$ by $w_q$.
Any two vertices connected by a path of edges marked by 
 $w_1,\dots,w_s$ correspond to generators conjugate 
 by the product $w_1w_2\dots w_s$.
Two vertices in $\Ga$ are in a common connected component
 if and only if the corresponding generators of $\Pi$ are conjugate.
\smallskip

By Definition~4.2 any Gauss diagram provides a Wirtinger presentation,
 where all generators split into classes of conjugate ones such that
 within each class the number of generators is equal 
 to the number of relators.
Some relators may follow from the remaining ones,
 as in the case of a classical knot.
After removing superfluous relators each connected component
 of $\Ga$ has Euler characteristic either 0 or 1, i.e.
 either it is a tree or it has exactly one cycle.
\end{proof}


\subsection{From a Wirtinger group to a virtual link}
\noindent
\smallskip

Here we prove the conditions of Proposition~4.5 are sufficient
 for the realizability of Gauss paragraphs. 
In Theorem~4.8 we construct a virtual link starting with 
 a suitable Wirtinger presentation.
The method is similar to \cite[section~3]{Kim}.
The first reduction of Lemma~4.6 converts a given presentation into a cyclic form, 
 where each relator says that two successive generators are conjugate.
The second reduction of Lemma~4.7 simplifies all conjugating elements
 and transforms them to generators.
\medskip

\noindent
{\bf Lemma~4.6.}
Let $\Pi$ be a Wirtinger presentation such that
 each connected component of the associated graph $\Ga$
 has Euler characteristic 0 or 1.
Then $\Pi$ can be converted into a presentation
 $\ab{m_1,\dots,m_n\vl r_1,\dots,r_n}$, where
 each relator has the form $r_q=m_i^{-1} w_q^{-1} m_{i+1} w_q$,
 $m_{n+1}=m_1$.
\begin{proof}
Let $\Ga'$ be a connected component of $\Ga$.
If $\Ga'$ has Euler characteristic 1 then add a superfluous relation 
 to the presentation to get a cycle (possibly, a loop) in $\Ga'$.
We will convert $\Ga'$ into a simple cycle, 
 where all vertices have degree~2.
After renumbering generators the resulting presentation
 will have a required cyclic form.
\smallskip

Suppose that $\Ga'$ has 2 edges $e_{ij},e_{jk}$ connecting
 the vertex $i$ with $j,k$ in such a way that 
 $e_{ij}$ is included into a cycle of $\Ga'$, 
 but $e_{jk}$ is not, see Fig.~12. 
Let the corresponding relations of $\Pi$ be 
 $m_i=w_q^{-1}m_jw_q$ and $m_j=w_p^{-1}m_kw_p$.
Then the first relation and edge $e_{ij}$ can be replaced
 by the relation $m_i=(w_pw_q)^{-1}m_k(w_pw_q)$ and 
 a new edge $e_{ik}\subset\Ga'$.
The procedure works until $\Ga'$ is a simple cycle.
The second reduction converts a given presentation 
 into a simple form, where each relator says that 
 two successive generators are conjugate by another generator.
\end{proof}

\begin{figure}[!h]
\includegraphics[scale=1.0]{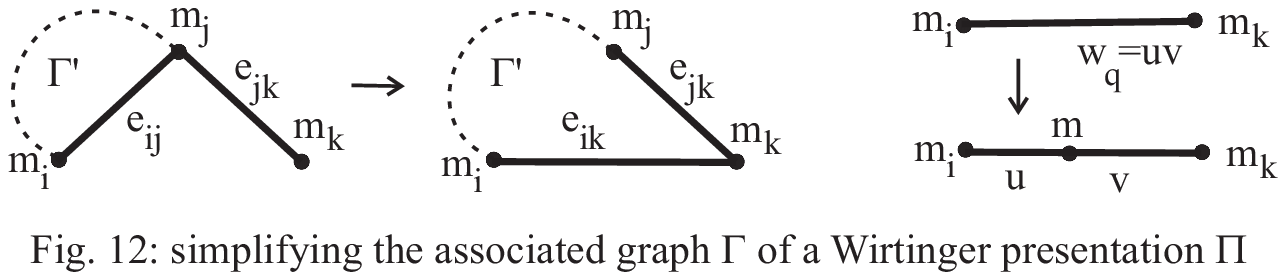}
\end{figure}

\noindent
{\bf Lemma~4.7.}
Let $\Pi$ be a Wirtinger presenation such that
 each connected component of the associated graph
 has Euler characteristic 0 or 1.
Then $\Pi$ can be transformed to a presentation
 $\ab{m_1,\dots,m_n\vl r_1,\dots,r_n}$, where
 each relator has the form $r_q=m_i^{-1} m_q^{-\e} m_{i+1} m_q^{\e}$,
 $\e=\pm 1$.
\begin{proof}
We introduce more generators and relators to reduce 
 the words $w_q$ from the given presentation to generators.
Assume that the conjugating word $w_q$ can be decomposed as 
 $w_q=uv$ for non-empty words $u,v$.
Then we add a new generator $m$ and replace the given relation 
 $m_i=w_q^{-1}m_jw_q$ by two new ones:
 $m_i=v^{-1}mv$ and $m=u^{-1}m_ju$.
\end{proof}
\smallskip

\noindent
{\bf Theorem~4.8.}
A group can be realized as a virtual link group
 if and only if it has a Wirtinger presentation 
 such that each connected component of 
 the associated graph has Euler characteristic 0 or 1.
\begin{proof}
The necessity was proved in Proposition~4.5.
Due to Lemmas~4.6--4.7 it suffices to construct
 a virtual link whose group has a given presentation
 $\Pi=\ab{m_1,\dots,m_n\vl r_1,\dots,r_n}$, where
 each relator has the form $r_q=m_i^{-1} m_q^{-1} m_{i+1} m_q$
 for some generator $m_q$, $q\in\{1,\dots,n\}$.
\smallskip

For each class of conjugate generators take a circle
 and split it into the same number of arcs.
Mark the arcs counterclockwisely by the successive generators.
For each relator $r_q=m_i^{-1} m_q^{-\e} m_{i+1} m_q^{\e}$,
 draw the arrow from a point on the arc marked by $m_q$
 to the common point of $m_i,m_{i+1}$.
The sign of the arrow is $\e$.
By Definition~4.2 the group of this Gauss diagram,
 i.e. of the required virtual link, has the given presentation $\Pi$.
\end{proof}
\smallskip

The Gauss paragraph $\{12^+,1^-2\}$ in Fig.~6 
 is not planar since its Carter surface is a torus.
The group of the associated Gauss diagram is isomorphic to $\Z^2$
 and is realizable as the group of the classical Hopf link,
 the graph $\Ga$ in Fig.~11 satisfies Theorem~4.8.
\smallskip

The criterion of Theorem~4.8 is similar to 
 the characterization of 2-dimensional link groups
 via their Wirtinger presenations \cite{Kam}. 
It is known that any abelian group can be realized
 as the second homology group of the fundamental group
 of the complement to a knotted surface in 4-sphere \cite{Lit}.
The same realizability question for virtual link groups
 is completely open and striking.
There is only one example of a virtual knot such that
 the second homology of its group is finite, $\Z_2$ \cite{KL}.



\end{document}